\documentclass[11pt]{article}

\usepackage{amsfonts}
\usepackage{amsmath}
\usepackage{amssymb}
\usepackage{latexsym,color}
\usepackage{amsbsy}
\usepackage{graphicx}

\newtheorem{thm}{Theorem}[section]
\newtheorem{con}{Conjecture}[section]
\newtheorem{lemma}[thm]{Lemma}
\newtheorem{cor}[thm]{Corollary}
\newtheorem{pro}[thm]{Proposition}
\newtheorem{example}[thm]{Example}
\newtheorem{definition}[thm]{Definition}

\newtheorem{remark}[thm]{Remark}
\newtheorem{Algorithm}[thm]{Algorithm}

\newcommand{\comment}[1]{}

\newcommand{\C}{{\mathbb C}}

\newcommand{\ncom}{\newcommand}
\ncom{\ns}{\normalsize}
\ncom{\la}{\lambda}
\ncom{\bm}{\boldmath}
\ncom{\noi}{\noindent}
\ncom{\bq}{\begin{equation}}
\ncom{\eq}{\end{equation}}
\ncom{\beqn}{\begin{eqnarray*}}
\ncom{\eeqn}{\end{eqnarray*}}
\ncom{\ba}{\begin{array}}

\ncom{\ea}{\end{array}}
\ncom{\beq}{\begin{eqnarray}}
\ncom{\eeq}{\end{eqnarray}}
\ncom{\nno}{\nonumber}
\ncom{\hs}{\mbox{\hspace{.25cm}}}
\ncom{\rar}{\rightarrow}
\ncom{\Rar}{\Rightarrow}
\ncom{\noin}{\noindent}
\ncom{\bc}{\begin{center}}
\ncom{\ec}{\end{center}}
\ncom{\sz}{\scriptsize}
\ncom{\fpd}{\Phi(\pi^{'})}
\ncom{\fp}{\Phi(\pi) }
\ncom{\nk}{\left< \begin{array}{c}
                       n\\k \end{array} \right>}
\ncom{\nd}{1^{'},2^{'},\cdots,n^{'}}
\ncom{\R}{I\!\!R}
\ncom{\de}{\bigtriangleup (F_{2n},\leq)}
\ncom{\del}{\bigtriangleup}
\ncom{\cov}{<\!\!\!\!\cdot }
\ncom{\bt}{\begin{thm}}
\ncom{\bcon}{\begin{con}}
\ncom{\et}{\end{thm}}
\ncom{\econ}{\end{con}}
\ncom{\bl}{\begin{lemma}}
\ncom{\el}{\end{lemma}}
\ncom{\bco}{\begin{cor}}
\ncom{\ds}{\displaystyle}
\ncom{\eco}{\end{cor}}
\ncom{\bp}{\begin{pro}}
\ncom{\ep}{\end{pro}}
\ncom{\bex}{\begin{example}}
\ncom{\eex}{\end{example}}
\ncom{\bd}{\begin{definition}}
\ncom{\ed}{\end{definition}}
\ncom{\brm}{\begin{remark}}
\ncom{\erm}{\end{remark}}
\ncom{\bal}{\begin{Algorithm}}
\ncom{\eal}{\end{Algorithm}}
\ncom{\ol}{\overline}

\ncom{\pf}{\noi {\bf Proof  }}
\ncom{\be}{\begin{enumerate}}
\ncom{\ee}{\end{enumerate}}
\ncom{\s}{\subset}
\ncom{\T}{{\cal T}}
\ncom{\B}{{\cal B}}
\ncom{\A}{{\cal A}}

\title{\Large{\textcolor{black} {\bf Note on Doron Zeilberger's
paper \\{\em${5\choose 2}$ proofs that ${n\choose k}\leq {n\choose {k+1}}$ 
if $k < n/2$
}}}}

\author{{\textcolor{black} {\bf Murali K. Srinivasan}} \\
{\em  \normalsize{Department of Mathematics}}\\
{\em  \normalsize{Indian Institute of Technology, Bombay}}\\
{\em  \normalsize{Powai, Mumbai 400076, INDIA}}\\
{\bf  \texttt{mks@math.iitb.ac.in,murali.k.srinivasan@gmail.com}}}

\begin{document}
\date{}
\maketitle
\begin{center}{\em To the memory of Mobi}\end{center}

In his very interesting paper {\bf \cite{z}}, Doron Zeilberger recalls,
among other results, several
proofs of unimodality of the binomial coefficients. Among the
combinatorial proofs he calls the beautiful symmetric chain decomposition (SCD) proof
of de Bruijn, Tengbergen, and Kruyswijk {\bf \cite{a,btk}} his personal
favourite. He then goes on to present a variation of Proctor's {\bf\cite{p}} 
beautiful algebraic proof (based on the representation theory of the Lie
algebra 
${\mbox{sl}}(2,\C)$) of injectivity of the up operator   
on the lower half of the Boolean algebra (which implies unimodality),
calling this the longest and yet the best proof.

The purpose of this note is to present yet another algebraic 
proof of the statement in
the title: we show  that the SCD proof has a simple and
natural linear analog that proves injectivity of the up
operator on the lower half (and surjectivity on the upper half).

For a finite set $S$, let $V(S)$ denote the complex vector space with $S$ as
basis. 
Let $B(n)$ denote the set of all subsets of $[n]=\{1,2,\ldots ,n\}$ and, for
$0\leq k \leq n$, let $B(n)_k$ denote the set of all $k$-subsets of $[n]$. 
Then we have
$ V(B(n))=V(B(n)_0)\oplus V(B(n)_1) \oplus \cdots \oplus V(B(n)_n)$ (vector space direct
sum).
An element $v\in V(B(n))$ is {\em homogeneous} if $v\in V(B(n)_k)$ for some
$k$. We say that a homogeneous element $v$ is of {\em rank}
$k$, and we write $r(v)=k$, if $v\in V(B(n)_k)$. 
%and we extend the notion of rank to homogeneous elements by writing $r(v)=i$.
%A linear map $T:V(P)\rar V(P)$ is said to be {\em order raising} if, for all
%$p\in P$, $T(p)$ is a linear combination of the elements covering $p$ (note
%that this implies that $T(p)=0$ for all maximal elements of $P$ and that
%$T(v)$ is homogeneous for homogeneous $v$). 
The {\em
up operator}  $U:V(B(n))\rar V(B(n))$ is defined, for $X\in B(n)$, by
$U(X)= \sum_{Y} Y$,
where the sum is over all $Y$ covering $X$, i.e., $X\subseteq Y$ and
$|Y|=|X|+1$.
A {\em symmetric Jordan chain} (SJC) in $V(B(n))$ is a sequence  
$v=(v_1,\ldots ,v_h)$ of nonzero homogeneous elements of $V(B(n))$
such that $U(v_{i-1})=v_i$, for
$i=2,\ldots h$, $U(v_h)=0$, and
$r(v_1) + r(v_h) = n$, if $h\geq
2$, or else $2r(v_1)= n$, if $h=1$. 
Note that the
elements of the sequence $v$ are linearly independent, being nonzero and of
different ranks. 
%We say that $v$ {\em
%starts} at rank $r(v_1)$ and {\em ends} at rank $r(v_h)$.
%If, in addition, $v$ is symmetric, i.e., 
%we say that $v$ is a
%{\em symmetric Jordan chain}.
A {\em symmetric Jordan basis} (SJB)  of $V(B(n))$ is a basis of $V(B(n))$
consisting of a disjoint union of SJC's in $V(B(n))$.

\noi
{\bf Theorem} {\em{There exists an SJB of $V(B(n))$}}.

\pf
We give a constructive proof that inductively produces an explicit SJB of
$V(B(n))$, the case $n=0$ being clear.

Let $V=V(B(n+1))$. Define $V(0)$ to be the subspace of $V$ generated by all
subsets of $[n+1]$ not containing $n+1$ and define $V(1)$ to be the subspace
of $V$ generated by all subsets of $[n+1]$ containing $n+1$. We have
$V=V(0)\oplus V(1)$. The linear map $R:V(0)\rar
V(1)$, given by  $X\mapsto X\cup \{n+1\},\;X\subseteq [n]$ is an isomorphism.
We write $R(v) = \ol{v}$. Let $U$ and $U_0$ denote, respectively, the up
operators on $V$ and $V(0)$ ($=V(B(n))$). We have, for $v\in V(0)$, 
\beq \label{a} 
&U(v)=U_0(v)+\ol{v},\;\;U(\ol{v})=\ol{U_0(v)}.
\eeq

By induction hypothesis there is an SJB $\cal{B}$ of $V(B(n))=V(0)$. We
shall now produce an SJB $\cal{B'}$ of $V$ by producing, for each
SJC in $\cal{B}$, either one or two SJC's in $V$ such that the
collection of all these SJC's is a basis.

Consider an SJC $(x_k,\ldots ,x_{n-k})$ (for some
$0\leq k \leq \lfloor n/2 \rfloor$),
where $r(x_k) = k$, in $\cal{B}$. 
%\beq \label{ba1}
%&(x_k,\ldots ,x_{n-k}),\;\;0\leq k \leq \lfloor n/2 \rfloor,&
%\eeq
%

%
We now consider two cases.

(a) $k=n-k$ : From (\ref{a}) we have $U(x_k)=\ol{x_k}$ and
$U(\ol{x_k})=\ol{U_0(x_k)}=0$. Since $R$ is an isomorphism $\ol{x_k}\not=0$. 
Add to $\cal{B'}$ the SJC
\beq \label{ba2}
&(x_k, \ol{x_k}).&
\eeq 

(b) $k< n-k$ : Set $x_{k-1}=x_{n+1-k}=0$ and define
%Add to $\cal{B'}$ the following two SJC's 
\beq \label{ba3}
&(y_k,\ldots ,y_{n+1-k}), \mbox{ and }
(z_{k+1},\ldots ,z_{n-k}), \eeq
by
\beq \label{ba4} 
y_l &=& x_l + (l-k)\, \ol{x_{l-1}},\;\;k\leq l \leq n+1-k. \\\label{ba5}
z_l &=& (n-k-l+1)\,\ol{x_{l-1}} -x_l,\;\;k+1\leq l \leq n-k.\eeq

From (\ref{a}) we have 
\beq \label{b}
&U(\ol{x_l})=\ol{U_0(x_l)} = \ol{x_{l+1}},\;k\leq l\leq n-k &
\eeq 

It thus follows from (\ref{a}) and (\ref{b}) that, for $k\leq l < n+1-k$, we have
$$U(y_l)=U(x_l + (l-k)\ol{x_{l-1}}) = x_{l+1} + \ol{x_l} +
(l-k)\ol{x_l}=x_{l+1}+(l-k+1)\ol{x_l}=y_{l+1}.$$ 
Note that when $l=k$ the second step above is justified 
because of the presence of the $(l-k)$ factor even though
$U(\ol{x_{k-1}})=0\not= \ol{x_k}$. We
also have $U(y_{n+1-k})=U((n+1)\ol{x_{n-k}})=(n+1)\ol{U_0(x_{n-k})}=0$.

Similarly, for $k+1\leq l < n-k$, we have
$$U(z_l)=U((n-k-l+1)\ol{x_{l-1}}-x_l) = (n-k-l+1)\ol{x_l} -x_{l+1} -
\ol{x_l} = (n-k-l)\ol{x_l} - x_{l+1} = z_{l+1}.$$ 
and $U(z_{n-k})=U(\ol{x_{n-k-1}} - x_{n-k})=\ol{x_{n-k}}-\ol{x_{n-k}}=0.$

Since $y_k = x_k\not= 0$, $y_{n+1-k} = (n+1)\ol{x_{n-k}}\not= 0$,
$x_l$ and $\ol{x_{l-1}}$ are linearly independent, for $k+1\leq l \leq n-k$
and the $2\times 2$ matrix 
$$\left(\ba{rl} 1 & l-k \\
                           -1 & n-k-l+1 \ea\right)$$
is nonsingular for $k+1\leq l \leq n-k$, it follows that (\ref{ba3}) gives
two independent SJC's in $V$. Add these two to $\cal{B'}$. 

Since $V=V(0)\oplus V(1)$ and $R$ is an isomorphism it follows that performing the above step for
each SJC in $\cal{B}$ we get an SJB $\cal{B'}$ of $V$. $\Box$

The main idea in the proof of the Theorem above  has several consequences,
studied in {\bf\cite{sr}}. In that paper we show the following:
\begin{itemize}
\item Introduce the standard inner product on $V(B(n))$, i.e., the
set $\{X\;|\;X\in B(n)\}$ is an orthonormal basis. It is shown
that the SJB produced above is orthogonal. Moreover, any two SJC's starting
at rank $k$ and ending at rank $n-k$ look alike in the sense that the ratios
of the lengths of the successive vectors is the same in both the chains and
these ratios (i.e., the singular values) can be explicitly written down. This yields a new constructive
proof of the explicit block diagonalization of the Terwilliger algebra of
the binary Hamming scheme, recently achieved by Schrijver {\bf\cite{sc}}.

\item It is shown that the SJB produced above is the canonically defined
(upto common scalars on each SJC) symmetric Gelfand-Tsetlin basis of $V(B(n))$. This
gives a natural representation-theoretic explanation for the orthogonality
in the item above.

\item  The algorithm given above can be generalized to produce an SJB in the
multiset case. The author's believes that the multiset version, which has a
deeper level of recursion than the set case, deserves further study from a
representation-theoretic viewpoint.

\end{itemize}

\begin{center} {\bf \Large{Acknowledgement}}
\end{center}
The author is grateful to Professor Doron Zeilberger whose encouraging
e-mail motivated him to write this note.

\end{document}